\documentclass{article}
\usepackage{amsfonts, amsmath,amssymb}
\newtheorem{thm}{Theorem}

\newcommand{\Fq}{{\mathbb F}_q}
\newcommand{\Fp}{{\mathbb F}_p}
\newcommand{\Fqast}{{\mathbb F}_q^\ast}
\def\deg{\mathop{\rm deg}}

\begin{document}

\title{Giesbrecht's algorithm, the HFE cryptosystem\\ and Ore's 
\lowercase{$p^s$}-polynomials}

\author{Robert S. Coulter\thanks{author now at Department of Mathematical Sciences, 
          University of Delaware, Newark, DE, 19716, 
          U.S.A.
A version of this paper appeared in the Proceedings of the 5th ASCM, 2001, pages 36--45, but is not widely available.
     }
\and George Havas \and Marie Henderson}  

\date{}
\maketitle

\begin{abstract}
We report on a recent implementation of Giesbrecht's algorithm for
factoring polynomials in a skew-polynomial ring. We also discuss the
equivalence between factoring polynomials in a skew-polynomial ring and
decomposing $p^s$-polynomials over a finite field, and how Giesbrecht's
algorithm is outlined in some detail by Ore in the 1930's. We end with some
observations on the security of the Hidden Field Equation
(HFE) cryptosystem, where $p$-polynomials play a central role.
\end{abstract}

\section{Introduction and Background}

Let $\Fq$ denote the finite field with $q=p^e$ elements, $p$ a prime. We
use $\Fqast$ to denote the non-zero elements of $\Fq$. The polynomial ring
in an indeterminate $X$ over any field $K$ will be denoted by $K[X]$ and
for $f,g \in K[X]$, $f\circ g=f(g)$ represents the composition of $f$ with
$g$. We recall that a {\em permutation} polynomial is a polynomial which
permutes the elements of the finite field under evaluation. A
$p$-{\em polynomial} (sometimes called an additive or linearised
polynomial) is a polynomial $L\in\Fq[X]$ of the shape   
\begin{equation*}
L(X) = \sum_{i} a_i X^{p^i}
\end{equation*}
with $a_i\in\Fq$. More specifically, for any integer $s$, a
$p^s$-polynomial is a $p$-polynomial where $a_i=0$ whenever $i$ is
not a multiple of $s$. We note that $p^s$-polynomials are closed under
composition (this is simply established).   

The problem of completely decomposing a polynomial $f \in K[X]$ into
indecomposable factors, where $K$ is a field, has a long and rich history.
When $K$ is the complex plane, Ritt~\cite{ritt22} showed that there
exists an essentially unique decomposition for any chosen
polynomial. It is unique in the sense that for any $f\in K[X]$
in a complete decomposition of $f$: the number of factors is invariant;
and the degrees of the factors are unique up to permutation.
So, if we have two complete decompositions 
\begin{align*}
f&=f_1\circ\cdots\circ f_m\\
&=g_1\circ\cdots\circ g_n,
\end{align*}
then $m=n$ and $\deg(f_i)=\deg(g_{\pi(i)})$ for some permutation $\pi$
of $\{1,\ldots,m\}$. Any class of polynomials defined over a field for
which this property holds is commonly said to satisfy Ritt's theorem. The
generalisation of Ritt's theorem to all fields of characteristic zero was
carried out by Engstrom\cite{eng41}, and Levi\cite{lev42}. However, for
fields of non-zero characteristic, the situation is not so clearcut.

A polynomial is called {\em wild} if its degree is divisible by the
characteristic $p$, and {\em tame} otherwise.
Any non-linear $p^s$-polynomial is therefore a wild polynomial.
A distinction between the behaviour of wild and tame polynomials arises
when one considers Ritt's theorem in the context of a finite
field. Fried and MacRae\cite{fried69} showed that any tame polynomial
satisfies Ritt's theorem. However, Dorey and Whaples\cite{dorey74}
gave an example which showed that not all wild polynomials satisfied
Ritt's Theorem. Other properties (not discussed in this article) of
tame and wild polynomials are also distinct. However, not
all wild polynomials deviate from tame polynomial
behaviour. Specific to this question, Ore\cite{ore33a} showed in the
1930's that $p$-polynomials satisfy Ritt's theorem. 

It is interesting to note that $p$-polynomials over a finite field
appear to be the second class of polynomials shown to satisfy
Ritt's theorem, after Ritt had established the complex field
case. This was not noted by Ore but is evident from his work: see
Ore\cite{ore33a} (Chapter 2, Theorem 4) which gives a statement
equivalent to Ritt's theorem for $p$-polynomials. A further class of
wild polynomials, known as $(p^s,d)$-polynomials (or, sub-linearised
polynomials) can be shown to satisfy Ritt's theorem by using results
of Henderson and Matthews\cite{henderson99}.  

Exponential-time algorithms for determining the complete decomposition
of polynomials were first given by Alagar and Thanh\cite{ala85}, and
Barton and Zippel\cite{bar85}. The first polynomial-time algorithm was
published by Kozen and Landau\cite{kozen89}, and separately by
Gutierrez, Recio and Ruiz de Velasco\cite{gut89}. These results were
improved for the tame case over a finite field by von zur
Gathen\cite{gathen90a}. A general purpose polynomial-time algorithm
for finding a complete decomposition of a rational function over an
arbitrary field was given by Zippel\cite{zippel91}. This last
algorithm provides a method for decomposing any polynomial, wild or
tame, over a finite field. However, one should note that in the wild
case, the algorithm simply finds any complete decomposition, as there
does not necessarily exist an essentially unique decomposition. 

\newpage
Although $p$-polynomials were the first polynomials over a finite
field shown to satisfy Ritt's theorem, they are the latest class of
polynomials for which a polynomial-time decomposition algorithm has
been given. The algorithm we refer to was described and analysed
by Giesbrecht\cite{giesbrecht98}. Giesbrecht presents his algorithm
in terms of factoring in skew-polynomial rings but it is well known
(and we later show) that the problem he considers is equivalent
to decomposing $p$-polynomials over a finite field. We
note that any decomposition algorithm for $p^s$-polynomials can
be adapted, at no computational cost, to decomposing
$(p^s,d)$-polynomials. For $(p,d)$-polynomials this was shown by
the authors\cite{coulter98b}, following earlier work
of Henderson and Matthews\cite{henderson99}. This can be extended to all
$(p^s,d)$-polynomials using the work of Ore\cite{ore33a}. This
subject is covered in another paper under preparation by the authors.

In this article, we report on a successful implementation of
Giesbrecht's algorithm, making some specific comments concerning the
probabilistic part of the algorithm. We also recall the work of
Oystein Ore, showing how Giesbrecht's algorithm is equivalent to an
algorithm described by Ore sixty years earlier. We also consider
implications of Ore's work to the security of the Hidden Field
Equations (HFE) cryptosystem. 

\section{Giesbrecht's algorithm and the work of Ore}

Giesbrecht\cite{giesbrecht98} introduces a probabilistic
polynomial-time algorithm for obtaining a complete (essentially
unique) factorisation of a polynomial in some classes of
skew-polynomial ring defined over a finite field. This problem is
intimately connected to the problem of determining an essentially
unique complete decomposition of $p$-polynomials, a class of {\em
wild} polynomials. In fact, there is a one-one correspondence between
factoring in a particular skew-polynomial ring over a finite field and
decomposing $p^s$-polynomials over a finite field. 

The skew-polynomial ring $\Fq[Y;\sigma]$, where $Y$ is an
indeterminate and $\sigma$ is an automorphism of $\Fq$, is a ring of
polynomials with the usual component-wise addition, and with
multiplication defined by $Ya =\sigma(a)Y$ for any $a\in\Fq$ (we
simply use juxtaposition to represent multiplication in $\Fq[X]$ and
$\Fq[Y;\sigma]$). Since $\sigma$ is an automorphism of $\Fq$, we must have
$\sigma(a)=a^{p^s}$ for some integer $s$. Given the definition of
multiplication above, it is easily seen that the skew-polynomial ring
$\Fq[Y;\sigma]$ is isomorphic to the ring of $p^s$-polynomials
over $\Fq$ with the operations of polynomial addition and composition.
Explicitly, the required isomorphism $\Phi$ satisfies $\Phi(X^p)\circ
\Phi(aX)=Ya=a^pY=\Phi(a^pX^p)$. From this it follows that
composition of $p^s$-polynomials acts in exactly the same manner as
multiplication in the skew-polynomial ring $\Fq[Y,\sigma]$. 

The theory introduced by Giesbrecht\cite{giesbrecht98} is developed in
its entirety in the works of Ore\cite{ore33a,ore33b,ore34}.
It may be more efficient to
implement Giesbrecht's algorithm using the $p^s$-polynomial
representation of the ring rather than the skew-polynomial ring
representation as set out in Giesbrecht's article but this is yet to
be tested. While Giesbrecht refers to Ore\cite{ore33b}, it is
in Ore's other two papers that he develops the algorithm
which Giesbrecht has rediscovered. Giesbrecht's key contribution is to
find a way of computing the crucial step, which is to find non-zero
zero divisors in a small algebra. He does this by using what he refers
to as Eigen rings. Ore\cite{ore33a} discusses the same method in
Chapter 2, Section 6 where he uses invariant rings. In particular,
Ore's Theorem 12 of that section is the key idea in Giesbrecht's
algorithm. Of course, Ore develops his theory in terms of
$p^s$-polynomials rather than skew-polynomial rings. Ore obtains these
results using an earlier paper, Ore\cite{ore33b}, where he developed
theory on factoring and primality of polynomials in more general
skew-polynomial rings than discussed here. The problem of
developing an algorithm for factoring polynomials over any
skew-polynomial ring remains open. 

Recently, a successful implementation of Giesbrecht's algorithm was
produced by Larissa Meinecke at the University of Queensland using the
Magma\cite{bos97} algebra package. There is one step in Giesbrecht's
algorithm which is probabilistic in nature, the rest of the algorithm
is strictly deterministic. Giesbrecht gives a lower bound for the
probability of this step being successful as 1/9. We have carried out
some testing regarding this step which suggests this lower bound is
very conservative. While we have been unable to determine a worst-case
scenario, in almost all cases tested, the step has been successful on
the first attempt.  

\section{HFE and $p$-polynomials}

The Hidden Field Equation (HFE) cryptosystem was introduced by
Patarin\cite{patarin96}. HFE is a public key cryptosystem and can be
described as follows: 
\begin{enumerate}
\item Choose a finite field $\Fq$, $q=p^e$, and a basis
$(\beta_1,\ldots,\beta_e)$ for $\Fq$ over $\Fp$. 
\item Select a polynomial $D$ of ``relatively small degree'' with the shape
\begin{equation*}
D(X) = \sum_{i,j} a_{ij} X^{p^i+p^j}
\end{equation*}
where $a_{ij} \in \Fq$ for all $i,j$. 
\item Choose two $p$-polynomials, $S$ and $T$, that permute $\Fq$.
\item Calculate $E(X) =S\circ D\circ T (X)\bmod (X^q - X)$.
\item Calculate $n_1,\ldots,n_e \in\Fp[X_1,\ldots,X_e]$ satisfying
\begin{equation*}
E(X) = \sum_{i=1}^e \beta_i n_i(X_1,\ldots,X_e)
\end{equation*}
and publish $\Fq$ and the $n_i$, $1 \le i \le e$. The polynomials $S, T$
and $D$ are the secret keys.
\end{enumerate}

If someone wishes to send a message $m$ to the owner of $E(X)$, then
they simply calculate $E(m)=y$ and send $y$. Decryption is carried out
by performing the following steps. As $S$ and $T$ are permutation
polynomials, they have functional and compositional (modulo $X^q-X$)
inverses. As $S$ and $T$ are known to the owner, they can determine
the inverse polynomials modulo $X^q-X$ (note that these inverses are also
$p$-polynomials). Thus the recipient of the message $y$ knows
$S$, $D$, $T$, $S^{-1}$ and $T^{-1}$. They determine $z$ satisfying
$S^{-1}(y)=z=D(T(m))$. Next they determine any $m_1\in\Fq$ so that
$D(m_1)=z$. Once $m_1$ is chosen they determine $m=T^{-1}(m_1)$. 
The middle step is only computationally feasible because the degree of
$D$ is chosen to be ``small''. 

The security of the system relies on
the assumption that if $\deg(E)$ is large, then solving for $m$ in
$E(m)=y$ is computationally infeasible. Note that several $m_1\in\Fq$
may need to be tried to find a ``sensible'' message $m$. This is
because $D$ is not necessarily chosen to be a permutation of $\Fq$ as
the authors of HFE assumed that this may be too difficult. However,
Blokhuis {\em et al.}\cite{coultera} have since given examples of
permutation polynomials from this class.  

Note that it makes no difference whether the polynomial $E$ or the set of $e$
polynomials $n_i$ is published if the basis used is known. In fact, an
attacker need not know the basis chosen as they may choose any basis 
to reconstruct a different, but effectively equivalent encryption function
(see the discussion below). If $E$ is constructed from the $e$ polynomials
$n_i$ using a different basis, alternative secret keys $S,T$ and $D$ may be
obtained and used to decipher messages.  

The HFE system is one of a family of cryptosystems which use
functional composition. Recently, some general attacks for these
systems were developed by Ye, Dai and Lam\cite{ye01}. An attack which
targets HFE specifically has been published by Kipnis and
Shamir\cite{kipnis99}. This is general in nature and is quite
successful, but does not break HFE in all cases. This attack has since
been improved by Courtois\cite{courtois01}. 

Polynomials with the shape $D$ are known as Dembowski-Ostrom
(DO) polynomials, see Dembowski\cite{dembowski68b}, Coulter and 
Matthews\cite{coulter97b} and Blokhuis {\em et al.}\cite{coultera}. 
For any $p$-polynomial $L\in\Fq[X]$ and any DO polynomial
$D\in\Fq[X]$, $L\circ D$ and $D\circ L$ are both DO polynomials. In
other words, DO polynomials are closed under composition with
$p$-polynomials. Also, it can be established that the reduction of a
DO polynomial modulo $X^q-X$ is again a DO polynomial. The HFE
description given above works in exactly the same way as that given by
Patarin\cite{patarin96} precisely because of the above comments,
coupled with the well known fact that any function over $\Fq$ can be
represented by a polynomial in $\Fq[X]$ of degree less than $q$ and a well
known result concerning linear operators (discussed below). 

Kipnis and Shamir\cite{kipnis99} note
several problems an attacker faces when they consider this scheme. We
address some of their concerns here. In the original description of HFE,
two linear transformations (or linear operators) over the vector space
${\mathbb F}_p^e$ are chosen, rather than two linearised polynomials as
described above. Kipnis and Shamir comment that ``these mixing operations
have natural interpretation over ${\mathbb F}_p$ but not over ${\mathbb
F}_{p^e}$, and it is not clear apriori that the $e$ published polynomials
over ${\mathbb F}_p$ can be described by a single univariate polynomial
$G$ over ${\mathbb F}_{p^e}$''. In fact, there is a natural
interpretation. Roman\cite{broman95} (pages 184-5) shows that every linear
operator on ${\mathbb F}_p^e$ can be represented by a linearised
polynomial over ${\mathbb F}_{p^e}$. So the description of HFE as given
above is equivalent. As DO polynomials are closed under composition with
linearised polynomials and their reduction modulo $X^q -X$ still results
in a DO, we are guaranteed that the published polynomials can be described
by a single univariate polynomial: it must be a DO. Kipnis and Shamir
continue ``Even if it exists (a single univariate polynomial), it may have
an exponential number of coefficients, and even if it is sparse, it may
have an exponentially large degree which makes it practically
unsolvable''. As the resulting polynomial is a DO polynomial, it has
$O(e^2)$ terms (compare to a random polynomial which has $O(p^e)$ terms),
which is not exponential. Certainly, the degree may be large. It remains
our objective, then, of finding a method of reducing the size of the
degree.

We can make more comments concerning the univariate description of HFE
given above. Let $E(X)$ be the public key, which is a DO
polynomial. Suppose we can determine $p$-polynomials $L_1$ and $L_2$
which are permutation polynomials and satisfy $L_1\circ f\circ L_2 =E$. 
Clearly, $f$ must also be a DO polynomial. Then we can
decrypt any message sent to the owner of $E$ using exactly the same
method used to decrypt in the standard way, but using the polynomials
$L_1, L_2$ and $f$, providing the degree of $f$ is sufficiently
small. Of course, it may not be possible to determine $p$-polynomials
that permute $\Fq$ which are left or right decompositional factors of
$E(X)$. However, when considering this problem, the following result
by Coulter and Matthews\cite{coulter97b}, immediately draws our
attention.  For any $a\in\Fq$ and any polynomial $t\in\Fq[X]$, define the
difference polynomial of $t$ with respect to $a$ by $\Delta_{t,a}(X) =
t(X+a) - t(X) - t(a)$. 

\begin{thm}\label{t1}
Let $f \in \Fq[X]$ with deg($f$) $<q$.
The following conditions are equivalent.
\par\noindent (i) $f = D + L$, where $D$ is a Dembowski-Ostrom
polynomial and $L$ is a $p$-polynomial.
\par\noindent (ii) For each $a\in\Fqast$, $\Delta_{f,a} = L_a$ where
$L_a$ is a $p$-polynomial depending on $a$. 
\end{thm}

\noindent This result provides an alternative definition of DO polynomials
and establishes an important connection between DO polynomials and
$p$-polynomials.

Let $E$ be the published DO polynomial used in the HFE
cryptosystem. We wish to find $L_1, L_2$ and $D$ satisfying
$E=L_1\circ D\circ L_2$. For the remainder, we set $f=D\circ L_2$ so
that $E=L_1\circ f$ and underline that $f$ is also a DO
polynomial. Our objective is to determine some information regarding
$L_1$. By Theorem~\ref{t1}, $\Delta_{E,a}$ is a
$p$-polynomial for any choice of $a$. Moreover, we have
\begin{align*}
\Delta_{E,a}(X) &= E(X+a) - E(X) - E(a)\\
&= L_1\left(f(X+a)\right) - L_1\left(f(X)\right) - L_1\left(f(a)\right)\\
&= L_1\left( f(X+a) - f(X) - f(a)\right)\\
&= L_1\circ \Delta_{f,a}.
\end{align*}
Thus for any non-zero choice of $a$, the polynomial $L_1$ is a left
decompositional factor of $\Delta_{E,a}$. Ore\cite{ore33b} shows that
there exists a left and right decomposition algorithm similar to the
well known greatest common divisor algorithm for a large class of
non-commutative polynomial rings (note that, in general, commutativity
for composition does not hold). He uses these results in Ore\cite{ore33a}
to establish and describe such algorithms for $p$-polynomials
specifically. In particular, using a variant of the Euclidean algorithm, we
can determine the Greatest Common Left-Decompositional Factor (GCLDF) of
two $p$-polynomials. This suggests the following method of attack to
determine the polynomial $L_1$.
\begin{enumerate}
\item Choose distinct elements $a_1,a_2\in\Fqast$.
\item Calculate $L(X)=$ GCLDF$(\Delta_{E,a_1}(X),\Delta_{E,a_2}(X))$.
\item Test to see if $L$ is a left decompositional factor of $E$. If it is,
then $L_1=L$ and we are done.
\item If $L$ is not a left decompositional factor of $E$, then choose a new
$a\in\Fqast$, distinct from previous choices, and calculate $L(X)=$
GCLDF$(L(X),\Delta_{E,a}(X))$. Return to Step 3. 
\end{enumerate}
We make the following observations. Step 3 can be carried out in time
$O(\log_p(\deg(E)))$ so has complexity much less than the Euclidean algorithm
calculation required in step 2 or 4. Note also that as Ore's work does not
extend to DO polynomials, one cannot simply calculate GCLDF$(L(X),E(X))$ to
obtain $L_1$. 

As mentioned, Giesbrecht's algorithm determines a complete decomposition of a
$p$-polynomial in probabilistic polynomial-time. However, this does not
mean we can determine $L_1$ methodically by completely decomposing $L$ after
step 2. Due to the nature of Ritt's theorem, we are not guaranteed that in a
full decomposition the proper factors of $L_1$ would be determined strictly on
the left. Further, the number of possible full decompositions is exponential in
the number of indecomposable factors. We make no claims at this point
concerning the number of GCLDF calculations required in step 4 to determine
$L_1$. It may require $O(q)$ such calculations, making the algorithm no better
than exhaustive search. Finally, we note that this attack does not
necessarily break HFE as the DO polynomial may not have a non-trivial GCLDF
and even if it did then the resulting DO polynomial may not be of
``sufficiently small'' degree. We are undertaking further research to
analyse this attack and to determine other methods of attacking HFE using
the connections between the DO polynomial and $p^s$-polynomial classes.
     
\section*{Acknowledgments}
This work was supported by the Australian Research Council.

\end{document}